\documentclass[11pt]{article}

\usepackage[letterpaper,left=1.25in,right=1.25in,top=1in,bottom=1in]{geometry}
\usepackage{amsmath,amssymb,amsthm}
\usepackage{setspace}
\usepackage{microtype}
\usepackage{needspace}
\usepackage{hyperref}

\allowdisplaybreaks
\newtheorem{theorem}{Theorem}[section]

\newtheorem{lemma}[theorem]{Lemma}
\newtheorem{claim}{Claim}[section]

\newtheorem{question}[theorem]{Question}
\newtheorem*{repeatedlowerbound}{Theorem~\ref{thm:lower}}

\newcommand{\tw}{\operatorname{tw}}
\newcommand{\ta}{\operatorname{tree}\text{-}\alpha}
\newcommand{\tc}{\operatorname{tree}\text{-}\chi}

\hypersetup{
  hidelinks,
  pdfauthor={Xinheng Lin},
  pdftitle={Treewidth and the Product of Tree-Independence and Tree-Chromatic Numbers},
  pdfsubject={Treewidth, tree-independence number, and tree-chromatic number},
  pdfkeywords={treewidth, tree-independence number, tree-chromatic number}
}

\title{Treewidth and the Product of
  Tree-Independence and Tree-Chromatic Numbers}
\author{Xinheng Lin\\[1.4ex]
  \normalsize Center for Discrete Mathematics, Fuzhou University\\[-0.1ex]
  \normalsize Fuzhou, P.\ R.\ China}
\date{October 13, 2026}

\begin{document}

\setlength{\abovedisplayskip}{6pt plus 2pt minus 2pt}
\setlength{\belowdisplayskip}{6pt plus 2pt minus 2pt}
\setlength{\abovedisplayshortskip}{0pt plus 2pt}
\setlength{\belowdisplayshortskip}{3pt plus 2pt minus 1pt}

\maketitle

\begin{abstract}
For a graph $G$, let $\tw(G)$, $\ta(G)$, and $\tc(G)$ denote its
treewidth, tree-independence number, and tree-chromatic number,
respectively.
We prove that every nonempty graph $G$ satisfies
\[
  \tw(G)+1\leq
  \ta(G)\tc(G)\bigl(1+\lfloor\log_2\ta(G)\rfloor\bigr).
\]
This improves the quadratic bound of Koutsoutis, Krause, Liu, Redzic,
and Ueckerdt~\cite{koutsoutis-et-al}. For every pair of integers
$n,k\geq2$, we construct a graph $G$ with $\ta(G)=n$, $\tc(G)=k$,
and $\tw(G)=\Omega(nk\log_2 n)$.
Thus the upper bound is best possible up to a constant factor.

\noindent\textbf{Key Words:} Treewidth, tree-independence number,
tree-chromatic number
\end{abstract}

\footnotetext[1]{2020 Mathematics Subject Classification: 05C15, 05C69, 05C75.\\
Email: \href{mailto:2500310007@fzu.edu.cn}{2500310007@fzu.edu.cn}}

\section{Introduction}

All graphs in this paper are finite and simple. For a graph $G$, let
$\alpha(G)$, $\omega(G)$, and $\chi(G)$ denote its independence number,
clique number, and chromatic number, respectively. For
$X\subseteq V(G)$, let $G[X]$ denote the subgraph of $G$ induced by $X$.

A \emph{tree-decomposition} of $G$ is a pair
$\mathcal T=(T,\{X_t\}_{t\in V(T)})$, where $T$ is a tree and
each $X_t\subseteq V(G)$ is a \emph{bag}, such that:
\begin{enumerate}
  \item for every edge $uv\in E(G)$, some bag $X_t$ contains both $u$ and $v$;
  \item for every vertex $v\in V(G)$, the set
  $\{t\in V(T):v\in X_t\}$ induces a nonempty subtree of $T$.\looseness=-1
\end{enumerate}
The \emph{width} of $\mathcal T$ is $\max_{t\in V(T)}|X_t|-1$.
The \emph{treewidth} of $G$, denoted by $\tw(G)$, is the minimum width
over all tree-decompositions of $G$.

Treewidth was introduced independently by
Halin~\cite{halin} in 1976 and by Robertson and
Seymour~\cite{robertson-seymour}, who developed it as a part of
the theory of graph minors.

The \emph{independence number} of a tree-decomposition
$\mathcal T=(T,\{X_t\}_{t\in V(T)})$ is
$\alpha(\mathcal T):=\max_{t\in V(T)}\alpha(G[X_t])$.
The minimum of $\alpha(\mathcal T)$ over all tree-decompositions
$\mathcal T$ of $G$ is its \emph{tree-independence number}, denoted
by $\ta(G)$, which was introduced independently by
Yolov~\cite{yolov} and by
Dallard, Milani\v{c}, and \v{S}torgel~\cite{dallard-milanic-storgel}.

Similarly, the \emph{chromatic number} of $\mathcal T$ is
$\chi(\mathcal T):=\max_{t\in V(T)}\chi(G[X_t])$.
The minimum of $\chi(\mathcal T)$ over all tree-decompositions
$\mathcal T$ of $G$ is its \emph{tree-chromatic number}, denoted
by $\tc(G)$, which was introduced by Seymour~\cite{seymour}.

Motivated by the elementary inequality $|V(G)|\leq\alpha(G)\chi(G)$,
Dallard, Milani\v{c}, and \v{S}torgel~\cite{dallard-milanic-storgel}
asked the following question.

\Needspace{6\baselineskip}
\begin{question}[Dallard, Milani\v{c}, and \v{S}torgel~{\cite[Question~8.4]{dallard-milanic-storgel}}]
\label{que:product}
Is the inequality
\[
  \tw(G)+1\leq\ta(G)\tc(G)
\]
valid for every graph $G$?
\end{question}

Koutsoutis, Krause, Liu, Redzic, and
Ueckerdt~\cite{koutsoutis-et-al} disproved the inequality in
Question~\ref{que:product} and established the following quadratic upper bound.

\Needspace{5\baselineskip}
\begin{theorem}[Koutsoutis et al.~{\cite[Proposition~3]{koutsoutis-et-al}}]
\label{thm:previous-upper}
For every graph $G$,
\[
  \tw(G)+1\leq\ta(G)^2\tc(G).
\]
\end{theorem}

Using Seymour's separation lemma~\cite[Theorem~1.1]{seymour},
we improve Theorem~\ref{thm:previous-upper} as follows.

\Needspace{5\baselineskip}
\begin{theorem}\label{thm:upper}
Every nonempty graph $G$ satisfies
\[
  \tw(G)+1\leq
  \ta(G)\tc(G)\bigl(1+\lfloor\log_2\ta(G)\rfloor\bigr).
\]
\end{theorem}

For the lower bound, Koutsoutis et al.~\cite{koutsoutis-et-al}
proved the following result.

\Needspace{9\baselineskip}
\begin{theorem}[Koutsoutis et al.~{\cite[Theorem~2]{koutsoutis-et-al}}]
\label{thm:previous-lower}
Let $f:\mathbb R\to\mathbb R_{>0}$ and
$h:\mathbb R_{>0}\to\mathbb R_{>0}$ be functions. Suppose that
$h(x)/x$ is nondecreasing and
$h(x)=o(x\log\log x)$ as $x\to\infty$.
Then some graph $G$ satisfies
\[
  \tw(G)>h\bigl(\ta(G)\bigr)f\bigl(\tc(G)\bigr).
\]
\end{theorem}

\Needspace{4\baselineskip}
Furthermore, we obtain the following lower bound.

\begin{theorem}\label{thm:lower}
There is an absolute constant $c>0$ such that, for all integers
$n,k\geq2$, some graph $G_{n,k}$ satisfies $\ta(G_{n,k})=n$,
$\tc(G_{n,k})=k$ and ${\tw(G_{n,k})+1\geq c\,nk\log_2 n}$.
\end{theorem}

Theorem~\ref{thm:lower} shows that the upper bound in
Theorem~\ref{thm:upper} is best possible up to a constant factor.

\Needspace{8\baselineskip}
\section{Proof of the upper bound}

A \emph{separation} of $G$ is a pair $(A,B)$ of subsets of $V(G)$ such
that $A\cup B=V(G)$ and there is no edge between $A\setminus B$ and
$B\setminus A$. Seymour~\cite[Theorem~1.1]{seymour}
proved the following lemma.

\begin{lemma}[Seymour~{\cite[Theorem~1.1]{seymour}}]\label{lem:seymour}
Every graph $G$ has a separation $(A,B)$ such that
$\chi(G[A\cap B])\leq\tc(G)$ and
\[
  \chi(G[A\setminus B]),\chi(G[B\setminus A])
  \geq\chi(G)-\tc(G).
\]
\end{lemma}

Using Lemma~\ref{lem:seymour}, we prove the following result.

\begin{lemma}\label{lem:coloring}
Every nonempty graph $H$ satisfies
\[
  \chi(H)\leq
  \tc(H)\bigl(1+\lfloor\log_2\alpha(H)\rfloor\bigr).
\]
\end{lemma}

\begin{proof}
We proceed by induction on $\alpha(H)$. If $\alpha(H)=1$, then
the conclusion is immediate.
We may therefore assume that $\alpha(H)\geq2$ and $\chi(H)>\tc(H)$.

By Lemma~\ref{lem:seymour}, there is a separation $(A,B)$ of $H$
such that both $H[A\setminus B]$ and $H[B\setminus A]$ have
chromatic number at least $\chi(H)-\tc(H)>0$ and hence are nonempty.
Since there are no edges between these two subgraphs, we have
$\alpha(H[A\setminus B])+\alpha(H[B\setminus A])\leq\alpha(H)$.
One of them, say $H'$, therefore satisfies
$1\leq\alpha(H')\leq\lfloor\alpha(H)/2\rfloor$. As
$\tc(H')\leq\tc(H)$, using Lemma~\ref{lem:seymour} again and the
induction hypothesis, we obtain
\[
\begin{aligned}
  \chi(H)&\leq\tc(H)+\chi(H')
   \leq\tc(H)+\tc(H')\bigl(1+\lfloor\log_2\alpha(H')\rfloor\bigr)\\
  &\leq\tc(H)
   \bigl(2+\lfloor\log_2\lfloor\alpha(H)/2\rfloor\rfloor\bigr)
   =\tc(H)\bigl(1+\lfloor\log_2\alpha(H)\rfloor\bigr),
\end{aligned}
\]
as required.
\end{proof}

\Needspace{4\baselineskip}
\begin{proof}[Proof of Theorem~\ref{thm:upper}]
Choose a tree-decomposition
$\mathcal T=(T,\{X_t\}_{t\in V(T)})$ of $G$ with
$\alpha(\mathcal T)=\ta(G)$. Let $X_t$ be a bag of maximum cardinality.
Then $\alpha(G[X_t])\leq\ta(G)$ and
$\tc(G[X_t])\leq\tc(G)$. By
Lemma~\ref{lem:coloring},
\begin{align*}
  |X_t|&\leq\alpha(G[X_t])\chi(G[X_t])
   \leq\alpha(G[X_t])\tc(G[X_t])
    \bigl(1+\lfloor\log_2\alpha(G[X_t])\rfloor\bigr)\\*
  &\leq\ta(G)\tc(G)\bigl(1+\lfloor\log_2\ta(G)\rfloor\bigr).
\end{align*}
By the definition of treewidth, we have
\[
  \tw(G)+1\leq
  \ta(G)\tc(G)\bigl(1+\lfloor\log_2\ta(G)\rfloor\bigr),
\]
as desired.
\end{proof}

\Needspace{8\baselineskip}
\section{A matching lower bound}\label{sec:lower}

The following construction can raise the treewidth to $|V(H)|-1$
while preserving the tree-chromatic number of the original graph $H$
and controlling the tree-independence number so that it does not go beyond
$\alpha(H)$.
We use the following definition and lemma from Koutsoutis, Krause,
Liu, Redzic, and Ueckerdt~\cite[Lemmas~5 and~6]{koutsoutis-et-al}.
For a graph $H$, its \emph{$1$-completion} $C(H)$ is obtained by
adding, for each pair of nonadjacent vertices $u,v\in V(H)$, a new
vertex whose neighborhood is exactly $\{u,v\}$.

\begin{lemma}[Koutsoutis et al.~{\cite[Lemmas~5 and~6]{koutsoutis-et-al}}]
\label{lem:completion}
Let $H$ be a graph with at least three vertices and at least one edge.
Then $\tc(C(H))=\tc(H)$, $\ta(C(H))\leq\alpha(H)$ and
$\tw(C(H))+1=|V(H)|$.
\end{lemma}

For a graph $G$ and a positive integer $q$, the \emph{lexicographic
product} $G[K_q]$, as defined in~\cite[Section~5.4]{hammack-imrich-klavzar},
is obtained by replacing each vertex $v\in V(G)$
by a clique $K_v$ of order $q$. For distinct $u,v\in V(G)$,
$K_u$ is complete to $K_v$ if and only if $uv\in E(G)$, and
anticomplete to $K_v$ otherwise.
For a positive integer $r$, the join $G\vee K_r$ is obtained by
adding a clique $K_r$ of order $r$ and making it complete to $G$;
for $r=0$, this means $G$ itself.

When $\tc(G)=\omega(G)$, we can increase the tree-chromatic number
and treewidth together while preserving the tree-independence number.
The following lemma extends the treewidth results in Theorem~24
(Section~4.5) and Theorem~22 (Section~4.2) of Gurski and
Weishaupt~\cite{gurski-weishaupt} to tree-independence and
tree-chromatic numbers.

\begin{lemma}\label{lem:clique-operations}
Let $G$ be a nonempty graph with $\tc(G)=\omega(G)$, and let
$q\geq1$ and $r\geq0$ be integers. If $H=G[K_q]\vee K_r$, then
$\ta(H)=\ta(G)$, $\tc(H)=q\tc(G)+r$ and
$\tw(H)+1=q\bigl(\tw(G)+1\bigr)+r$.
\end{lemma}

\begin{proof}
We first prove the first equality. To show $\ta(H)\leq\ta(G)$,
choose a tree-decomposition
$(T_G,\{U_t\}_{t\in V(T_G)})$ of $G$ with nonempty bags such that
$\alpha(G[U_t])\leq\ta(G)$ for every $t\in V(T_G)$.
Set $T_H:=T_G$.
For each $t\in V(T_G)$, put
$V_t:=V(K_r)\cup\bigcup_{v\in U_t}V(K_v)$.
Then $(T_H,\{V_t\}_{t\in V(T_H)})$ is a tree-decomposition of $H$.
Since $\alpha(G[U_t])=\alpha(H[V_t])$ for every $t\in V(T_G)$,
we have
\[
  \ta(H)\leq\max_{t\in V(T_H)}\alpha(H[V_t])
  =\max_{t\in V(T_G)}\alpha(G[U_t])=\ta(G).
\]

\Needspace{5\baselineskip}
Since $G$ is isomorphic to an induced subgraph of $H$, we have
$\ta(H)\geq\ta(G)$. Combining this with the preceding inequality
gives $\ta(H)=\ta(G)$.

\smallskip
We next prove the second equality. To show $\tc(H)\leq q\tc(G)+r$,
choose a tree-decomposition $(T'_G,\{W_t\}_{t\in V(T'_G)})$ of $G$
such that $\chi(G[W_t])\leq\tc(G)$ for every $t\in V(T'_G)$.
Set $T'_H:=T'_G$. For each $t\in V(T'_G)$, put
$X_t:=V(K_r)\cup\bigcup_{v\in W_t}V(K_v)$.
$(T'_H,\{X_t\}_{t\in V(T'_H)})$ is also a
tree-decomposition of $H$.
We have
\[
  \tc(H)\leq\max_{t\in V(T'_H)}\chi(H[X_t])
  \leq q\max_{t\in V(T'_G)}\chi(G[W_t])+r
  =q\tc(G)+r.
\]

To show $\tc(H)\geq q\tc(G)+r$, recall that every clique is
contained in a bag of every tree-decomposition (see
Diestel~\cite[Section~12.3, Corollary~12.3.5]{diestel} for a short proof).
Using the hypothesis $\tc(G)=\omega(G)$, we obtain
\[
  \tc(H)\geq\omega(H)=q\omega(G)+r=q\tc(G)+r.
\]
Combining the two inequalities gives $\tc(H)=q\tc(G)+r$.

\smallskip
Finally, we prove the last equality. To show
$\tw(H)+1\leq q(\tw(G)+1)+r$,
choose a tree-decomposition $(T''_G,\{Y_t\}_{t\in V(T''_G)})$ of $G$
of width $\tw(G)$. Set $T''_H:=T''_G$. For each $t\in V(T''_G)$, put
$Z_t:=V(K_r)\cup\bigcup_{v\in Y_t}V(K_v)$.
Again, $(T''_H,\{Z_t\}_{t\in V(T''_H)})$ is a
tree-decomposition of $H$. We have
\[
  \tw(H)+1\leq\max_{t\in V(T''_H)}|Z_t|
  =q\max_{t\in V(T''_G)}|Y_t|+r
  =q\bigl(\tw(G)+1\bigr)+r.
\]

To show $\tw(H)+1\geq q(\tw(G)+1)+r$, take a
tree-decomposition $(T'''_H,\{A_t\}_{t\in V(T'''_H)})$ of $H$ of width
$\tw(H)$. Let $T'''_G$ be the subgraph of $T'''_H$ induced by
$\{t\in V(T'''_H):V(K_r)\subseteq A_t\}$.
If $r=0$, then $T'''_G=T'''_H$. Using the fact that $K_r$
is a clique and the second condition of a tree-decomposition, we
conclude that $T'''_G$ is a nonempty subtree of $T'''_H$.
For each $t\in V(T'''_G)$, put
$B_t:=\{v\in V(G):V(K_v)\subseteq A_t\}$.
For each $v\in V(G)$, the nodes $t$ with $v\in B_t$ form a
nonempty intersection of subtrees, since $K_v\vee K_r$ is a clique.
For each $uv\in E(G)$, some bag $A_t$ contains the clique
$K_u\vee K_v\vee K_r$, so $u,v\in B_t$.
Thus $(T'''_G,\{B_t\}_{t\in V(T'''_G)})$ is a tree-decomposition of $G$.

Choose $t_0\in V(T'''_G)$ such that $|B_{t_0}|$ is maximum. By the
definition of treewidth,
\[
  \tw(G)+1\leq\max_{t\in V(T'''_G)}|B_t|=|B_{t_0}|.
\]
\begin{samepage}
The bag $A_{t_0}$ contains the pairwise disjoint vertex sets $V(K_v)$
for all $v\in B_{t_0}$ and $V(K_r)$.
Since each $K_v$ has order $q$, we have
$q|B_{t_0}|+r\leq|A_{t_0}|$. The chosen tree-decomposition of
$H$ has width $\tw(H)$, so $|A_{t_0}|\leq\tw(H)+1$.
Consequently,
\[
  q\bigl(\tw(G)+1\bigr)+r
  \leq q|B_{t_0}|+r
  \leq |A_{t_0}|
  \leq\tw(H)+1.
\]
Combining the two inequalities gives
$\tw(H)+1=q(\tw(G)+1)+r$.
\end{samepage}
\end{proof}

\Needspace{8\baselineskip}
The following construction is inspired by the shift graphs in Huynh,
Reed, Wood, and Yepremyan~\cite[Section~2]{huynh-reed-wood-yepremyan}.

For convenience, we restate Theorem~\ref{thm:lower} here.

\begin{repeatedlowerbound}
There is an absolute constant $c>0$ such that, for all integers
$n,k\geq2$, some graph $G_{n,k}$ satisfies $\ta(G_{n,k})=n$,
$\tc(G_{n,k})=k$ and ${\tw(G_{n,k})+1\geq c\,nk\log_2 n}$.
\end{repeatedlowerbound}

\Needspace{4\baselineskip}
\begin{proof}
\renewcommand{\theclaim}{\getrefnumber{thm:lower}.\arabic{claim}}
\setcounter{claim}{0}
Fix an integer $n\geq3$. Let
\[
  \mathcal I_n:=\bigl\{[x,x+3^i]:
  i,x\in\mathbb Z,\ 0\leq i\leq\lfloor\log_3 n\rfloor,\ 0\leq x\leq n-3^i\bigr\}.
\]
\Needspace{4\baselineskip}
\noindent
Define $H_n$ to be the graph with vertex set $\mathcal I_n$ in which two
distinct vertices are adjacent if and only if they can be written as
$[x,y]$ and $[z,w]$ such that
\begin{equation}\label{eq:adjacency}
  x<z\leq y<w
  \qquad\text{and}\qquad
  y-z\in\{0,1,3,9,\ldots\}.
\end{equation}
\begin{claim}\label{clm:independence}
$\alpha(H_n)\leq n$.
\end{claim}

\begin{proof}
Let $\mathcal S$ be an independent set of $H_n$. Form a graph $F$ on
$\{0,\ldots,n\}$ by including the edge $xy$ for every $[x,y]\in\mathcal S$,
and orient each edge from its smaller endpoint to its larger endpoint.
There is no directed path $x\to y\to z$ of length two, since
$[x,y]$ and $[y,z]$ would be adjacent in $H_n$. Hence every nonisolated
vertex of $F$ is either a source or a sink.

Suppose that $C$ is a cycle in $F$. Choose
the smallest sink $u$ on $C$, and let $\ell<r$ be its two neighbors on
$C$. Let $v$ be the other neighbor of $r$ on $C$. The vertex $r$ is a source,
so $v$ is a sink, and the choice of $u$ gives
$\ell<r<u<v$. Since $[r,u]$ is a vertex of $H_n$, the construction
of $H_n$ implies that $u-r$ is a power of three. By
\eqref{eq:adjacency}, the vertices $[\ell,u]$ and $[r,v]$ are adjacent,
a contradiction. Thus $F$ is a forest, and
$|\mathcal S|=|E(F)|\leq n$.
\end{proof}

\Needspace{5\baselineskip}
\begin{claim}\label{clm:path}
$H_n$ has a tree-decomposition of chromatic number at most two.
\end{claim}

\begin{proof}
Let $P$ be the path with vertices $0,\ldots,n$ in this order. For each
$t\in V(P)$, put $X_t:=\{[x,y]\in V(H_n):x\leq t\leq y\}$.
Clearly,
$(P,\{X_t\}_{t\in V(P)})$ is a tree-decomposition of $H_n$.

We show that $\chi(H_n[X_t])\leq2$ for every $t\in V(P)$.
Fix $t\in V(P)$. Give a vertex $[x,y]\in X_t$ color $1$ if
$t<(x+y)/2$, and color $2$ if $t>(x+y)/2$.
Since $y-x$ is odd and $t$ is an integer, $(x+y)/2\neq t$, so
every vertex receives a color.

Let $[x,y]$ and $[z,w]$ be adjacent vertices in $X_t$, labelled so that
$x<z\leq y<w$. Since both vertices belong to $X_t$, we have
$z\leq t\leq y$. If $z=t=y$, then
$(x+y)/2<t<(z+w)/2$.
Otherwise, $z<y$. By \eqref{eq:adjacency} and the definition of $H_n$,
$y-x$, $y-z$, and $w-z$ are powers of three. $y-z$ is smaller
than both $y-x$ and $w-z$, so each interval has length at least three times
$y-z$, and hence $(x+y)/2<t<(z+w)/2$.
Thus $[x,y]$ receives color $2$ and $[z,w]$ receives color $1$.
Therefore, this coloring is proper.
Thus $\chi(H_n[X_t])\leq2$ for every $t\in V(P)$, as required.
\end{proof}

Counting intervals of lengths $3^i$ with
$0\leq i<\lfloor\log_3 n\rfloor$ gives
\[
  |V(H_n)|\geq\sum_{i=0}^{\lfloor\log_3 n\rfloor-1}(n+1-3^i)
    \geq\sum_{i=0}^{\lfloor\log_3 n\rfloor-1}\frac{2n}{3}
    =\frac{2n}{3}\lfloor\log_3 n\rfloor
    \geq\frac{n\log_2 n}{3\log_2 3},
\]
since $i<\lfloor\log_3 n\rfloor$ implies $3^i\leq n/3$, and
$\lfloor\log_3 n\rfloor\geq\frac12\log_3 n$ for $n\geq3$.

The vertices $[0,1]$ and $[1,2]$ are adjacent, so
Claim~\ref{clm:path} gives $\tc(H_n)=2$. Applying
Lemma~\ref{lem:completion} together with
Claim~\ref{clm:independence} yields $\tc(C(H_n))=2$,
$\ta(C(H_n))\leq n$ and $\tw(C(H_n))+1=|V(H_n)|$.

\Needspace{4\baselineskip}
For $n\geq3$, let $B_n$ be the disjoint union of $C(H_n)$ and $K_{n,n}$.
Since $\tc(K_{n,n})=2$, $\tw(K_{n,n})=n$ and
$\ta(K_{n,n})=n$, taking the maximum over components gives
$\ta(B_n)=n$, $\tc(B_n)=2$ and
$\tw(B_n)+1=\max\{|V(H_n)|,{n+1}\}$.
There are $n$ intervals of length one and $n-2$ of length three, so
$|V(H_n)|\geq2n-2\geq n+1$. Hence $\tw(B_n)+1=|V(H_n)|$.
Also, $B_n$ has an edge and $\omega(B_n)\leq\tc(B_n)=2$, so
$\omega(B_n)=2$.

For $n=2$, set $B_2:=K_{2,2}$, so $\ta(B_2)=\tc(B_2)=2$ and
$\tw(B_2)+1=3$. Together with the lower bound on $|V(H_n)|$, this gives,
for every integer $n\geq2$,
$\ta(B_n)=n$, $\tc(B_n)=\omega(B_n)=2$ and
$\tw(B_n)+1\geq n\log_2 n/(3\log_2 3)$.

Now fix integers $n\geq2$ and $k\geq2$, and put
$q:=\lfloor k/2\rfloor$ and $r:=k-2q\in\{0,1\}$. Define
$G_{n,k}:=B_n[K_q]\vee K_r$.
Since $\tc(B_n)=\omega(B_n)=2$, Lemma~\ref{lem:clique-operations}
gives $\ta(G_{n,k})=n$ and
$\tc(G_{n,k})=2q+r=k$.
\begin{samepage}
Since $q\geq k/3$ and $r\geq0$, the treewidth formula in the same lemma gives
\[
  \tw(G_{n,k})+1=q\bigl(\tw(B_n)+1\bigr)+r
   \geq\frac{k}{3}\cdot\frac{n\log_2 n}{3\log_2 3}
   =\frac{nk\log_2 n}{9\log_2 3}.
\]
This proves the theorem with $c=1/(9\log_2 3)$.
\end{samepage}
\end{proof}

\Needspace{6\baselineskip}
\medskip
\noindent\textbf{Statement of AI Use.}
The central proof idea in this manuscript was developed with the assistance of
GPT 6.1 Astra. The author subsequently refined the idea and wrote
the manuscript, using the same model to assist with language polishing.
The author takes full responsibility for the content of this paper.

\end{document}